\documentclass[12pt]{amsart}
\usepackage{amssymb,amscd,amsthm,verbatim}
\usepackage{amsbsy}
\usepackage{latexsym}
\usepackage[all,cmtip]{xy}
\usepackage[mathscr]{euscript}
\usepackage{bbm}
\usepackage{stmaryrd}
\usepackage{cite}
\usepackage{tikz}

\date{\today}

\newcommand{\C}{{\mathbbm{C}}}

\newcommand{\R}{{\mathbbm{R}}}
\newcommand{\T}{{\mathbbm{T}}}

\newcommand{\Z}{{\mathbbm{Z}}}

\newcommand{\cH}{{\mathcal{H}}}

\newcommand{\cS}{{\mathcal{S}}}

\newcommand{\LP}{{\mathrm{LP}}}

\newcommand{\tr}{{\mathrm{Tr}}}

\newtheorem{theorem}{Theorem}[section]

\newtheorem{conj}[theorem]{Conjecture}

\theoremstyle{definition}

\newtheorem{question}{Question}
\newtheorem{issue}{Fundamental Issue}

\theoremstyle{plain}

\allowdisplaybreaks \numberwithin{equation}{section}

\newcommand{\set}[1]{\left\{#1\right\}}

\begin{document}

\title{Schr\"odinger Operators with Thin Spectra}
\begin{abstract}
The determination of the spectrum of a Schr\"odinger operator is a fundamental problem in mathematical quantum mechanics.
We discuss a series of results showing that Schr\"odinger operators can exhibit spectra that are remarkably thin in the sense of Lebesgue measure and fractal dimensions.
We begin with a brief discussion of results in the periodic theory, and then move to a discussion of aperiodic models with thin spectra.
\end{abstract}

\author[D.\ Damanik]{David Damanik}
\address{Department of Mathematics, Rice University, Houston, TX~77005, USA}
\email{damanik@rice.edu}
\thanks{D.D.\ was supported in part by NSF grant DMS--1700131 and by an Alexander von Humboldt Foundation research award}

\author[J.\ Fillman]{Jake Fillman}
\address{Department of Mathematics, Texas State University, San Marcos, TX 78666, USA}
\email{fillman@txstate.edu}
\thanks{J.F.\ was supported in part by a Simons Collaboration Grant}

\maketitle
%

\section{Introduction}

The purpose of this text is to discuss some recent work on Schr\"odinger operators with thin spectra. We will begin with a list of motivating questions and then focus on one of them, perhaps the most fundamental among them, and progressively ask more and more detailed refinements of this question, each of them prompted by a discussion of the previous one. This progression will lead us to the current state of affairs, where Schr\"odinger operators with surprisingly thin spectra have been identified, likely somewhat at odds with the expectations one might have had based on intuition or classical paradigms.

This note is not meant to be a comprehensive survey, and the questions we ask and the results we present will be selected according to the goal outlined in the previous paragraph: motivating and explaining the recent results on Schr\"odinger operators with thin spectra. For additional references and a discussion of the historical context we refer the reader to the original articles. Throughout the text, where appropriate, we point the reader towards suitable survey articles and textbooks for additional reading.

\section{Schr\"odinger Operators and Things we Might Want to Know About Them}

Schr\"odinger operators are central to the mathematical formulation of quantum mechanics, indeed their associated unitary group generates the time evolution of a quantum mechanical system. We will focus on the simplest case of a single (quantum) particle exposed to a potential $V$ and ignore all physical constants. Doing so we arrive at operators of the form $H=H_V$,
\begin{equation}\label{e.schrodingeroperator}
H = -\Delta + V
\end{equation}
in
$$
\cH = L^2(\R^d),
$$
where $\Delta$ is the standard \emph{Laplacian},
\begin{equation}\label{e.laplacian}
[\Delta \psi](x) = \sum_{j = 1}^d \frac{\partial^2 \psi}{\partial x_j^2}(x),
\end{equation}
and $V$ acts by multiplication with the \emph{potential} $V : \R^d \to \R$ (denoted by the same symbol, as is customary),
\begin{equation}\label{e.potential}
[V \psi](x) = V(x) \psi(x).
\end{equation}

We will refer to $d$ as the \emph{dimension}. Its value is of crucial importance; there is a huge difference between the cases $d = 1$ and $d > 1$, both in terms of the results one should expect to hold and the methods that exist to establish these results. Even in the \emph{higher-dimensional case} $d > 1$, the value of $d$ is sometimes of importance, and there may be further transitions in the expected or provable behavior as $d$ is increased.

It is often beneficial to initially consider the \emph{free case} where $V$ vanishes identically, so that one deals with $-\Delta$, and to only later add on the potential, essentially as a perturbation of $-\Delta$.

Naturally one cannot define \eqref{e.laplacian} on all of $\cH$, at least if one wants that object to belong to $\cH$ as well. Thus, it is an important problem to identify a suitable \emph{domain} of this operator. It should not be too large (so that $-\Delta$ does not map some elements of the proposed domain outside $\cH$) and it should also not be too small (for reasons that become important later: we want to have a well defined time evolution $e^{-itH}$). This leads to the goal of choosing the domain such that the operator is \emph{self-adjoint}. Once this has been accomplished for the free case, the same domain works for large classes of potentials $V$, and hence the two-step procedure works out neatly in these cases. We won't dwell on this particular issue further, referring the reader instead to \cite{ReedSimon2}, and will assume henceforth that self-adjointness will always have been addressed successfully for the Schr\"odinger operators $H$ we consider. In a similar vein, most theorems may be formulated for various classes such as continuous potentials, smooth potentials, uniformly locally square-integrable potentials, and could even be pushed beyond that to even weaker regularity classes; to keep the exposition free of technicalities, we will mostly suppress assumptions on the regularity of $V$, and we will not attempt to state theorems in their absolute most general form. The reader is welcome to suppose that all potentials are bounded and continuous unless they are explicitly said to be otherwise.

Having realized \eqref{e.schrodingeroperator} as a self-adjoint operator in $\cH$, we can infer that the \emph{spectrum} $\sigma(H)$ must be real. We will comment on the importance of this set momentarily. Let us first note that several interesting decompositions of this set exist. The most basic decomposition is
$$
\sigma(H) = \sigma_\mathrm{disc}(H) \cup \sigma_\mathrm{ess}(H),
$$
where the \emph{discrete spectrum} $\sigma_\mathrm{disc}(H)$ consists of all eigenvalues of $H$ that have finite multiplicity, and the \emph{essential spectrum} is the complementary set. In particular, this is a disjoint union, that is, a partition of the spectrum. By definition, and the fact that $L^2(\R^d)$ is separable, the discrete spectrum is a countable set, and it is therefore irrelevant for the questions we will discuss here -- whether the spectrum is ``thin.'' Another important aspect, which however is less directly related to our theme, is that the dependence of $\sigma_\mathrm{disc}(H)$ on $V$ is very sensitive, while the dependence of $\sigma_\mathrm{ess}(H)$ on $V$ is far more robust.

In any event, let us formulate the first fundamental issue:

\begin{issue}
Study the essential spectrum $\cS = \sigma_\mathrm{ess}(H)$.
\end{issue}

Another consequence of self-adjointness is the existence of spectral measures guaranteed by the spectral theorem. Specifically, if $H$ is self-adjoint and $\psi \in \cH$, then there is a finite Borel measure $\mu_{H,\psi}$ on $\R$ (indeed, $\mu_{H,\psi} (\R) = \|\psi\|^2$) such that
\begin{equation}\label{e.specmeasure}
\langle \psi, (H-z)^{-1} \psi \rangle = \int_\R \frac{d\mu_{H,\psi}(E)}{E - z}, \quad z \not\in \sigma(H).
\end{equation}
In fact, each of the measures $\mu_{H,\psi}$ is supported by the spectrum, and hence we could integrate in \eqref{e.specmeasure} over $\sigma(H)$ rather than all of $\R$. For this reason, the structure of the spectrum places restrictions on the spectral measures.

For example, it is quite standard to decompose a Borel measure $\mu$ on $\R$ into its pure point part $\mu^\mathrm{(pp)}$ (which has a countable support), its singular continuous part $\mu^\mathrm{(sc)}$ (which gives no weight to countable sets and has a support of zero Lebesgue measure), and its absolutely continuous part $\mu^\mathrm{(ac)}$ (which gives no weight to sets of zero Lebesgue measure). Thus, our second fundamental issue is the following:

\begin{issue}
Study the spectral measures $\mu_{H,\psi}$ and their Lebesgue decomposition.
\end{issue}

As a specific example of how information regarding the first fundamental issue may give information regarding the second fundamental issue, imagine that $H$ is such that $\cS$ can be shown to have zero Lebesgue measure -- it then immediately follows that $\mu^\mathrm{(ac)}_{H,\psi}$ vanishes for every $\psi \in \cH$.

Why would the latter statement be of interest? As alluded to above, the time evolution of interest is given by the unitary group generated by $H$. In other words, since the time-dependent Schr\"odinger equation is given by
\begin{equation}\label{e.schrodingerequation}
i \frac{\partial}{\partial t} \psi = H \psi, \quad \psi(0) = \psi_0,
\end{equation}
self-adjointness of $H$ and the spectral theorem again ensure that we can write the solution as
\begin{equation}\label{e.schrodingerequationsolution}
\psi(t) = e^{-itH} \psi_0.
\end{equation}

The third fundamental issue is therefore the following:
\begin{issue}
Study the quantum evolution $e^{-itH} \psi_0$.
\end{issue}

There are several results (the RAGE Theorem~\cite{CFKS} being the most prominent) that link the second and third fundamental issue. The Lebesgue decomposition of $\mu^\mathrm{(ac)}_{H,\psi_0}$ (e.g., information about which of the three components are non-trivial) gives information about the behavior of $e^{-itH} \psi_0$ as $|t| \to \infty$. As a rather obvious instance of this connection, let us note that by \eqref{e.specmeasure} and \eqref{e.schrodingerequationsolution}, the so-called \emph{return probability}
$$
|\langle \psi_0, e^{-itH} \psi_0 \rangle|^2 = \Big| \int_{\sigma(H)} e^{-itE} \, d\mu_{H,\psi_0}(E) \Big|^2
$$
will go to zero as $|t| \to \infty$ if $\mu_{H,\psi_0} =\mu^\mathrm{(ac)}_{H,\psi_0}$ by the Riemann-Lebesgue Lemma.

To summarize, while we are ultimately interested in addressing the third fundamental issue, the standard approach to doing so involves addressing the other two fundamental issues as well. In this sense, a study of the (essential) spectrum of $H$ is the most basic of these tasks, and it is the one to which we will devote the remainder of this text.

\section{The Essential Spectrum $\cS$ and Things we Might Want to Know About it}

As explained in the previous section, we wish to study the essential spectrum $\cS = \sigma_\mathrm{ess}(H)$ of a Schr\"odinger operator $H$ in $\cH$. Inspired by the two-step procedure mentioned earlier, let us begin with the free case.

If $V$ vanishes identically, one can see via the Fourier transform, which conjugates $-\Delta$ to multiplication with $|\xi|^2$ in $L^2(\R^d, d\xi)$, that $\sigma(H) = \sigma_\mathrm{ess}(H) = \cS = [0,\infty)$ since that is the (closure of the) range of the multiplication operator in question. Thus, in the free case, $\cS$ is connected, that is, it has no gaps.

Adding on a potential, one may first wonder whether this property is preserved. If the operator $V$ is relatively compact with respect to $-\Delta$, which will hold if the function $V$ decays at infinity in a suitable sense, then the essential spectrum is indeed preserved and continues to be $[0,\infty)$.
If relative compactness fails, one cannot immediately deduce that the essential spectrum is not preserved, but one may confidently expect that this should at least be possible for some such perturbations of $-\Delta$. Thus, one should expect a positive answer to the following question:

\begin{question}\label{q.disconnected}
Can $\cS$ have gaps, that is, can it be disconnected?
\end{question}

We will describe a mechanism that leads to a positive answer in Section~\ref{s.disconnected}: in one dimension, it suffices to consider (non-constant) periodic potentials $V$. These potentials do not decay at infinity but they have a sufficiently simple structure that their essential spectra can be studied quite effectively.

Once it is understood that $\cS$ may have gaps, one may then ask if one can have arbitrarily many, or even infinitely many, gaps. Thus, we refine the previous question as follows:

\begin{question}\label{q.manygaps}
Can $\cS$ have infinitely many gaps/infinitely many connected components?
\end{question}

This question can be answered again in the framework of non-constant periodic potentials, but there is an additional wrinkle: the dimension $d$ now matters a lot! When $d = 1$, it is quite easy to produce infinitely many gaps, while in the case $d > 1$ it is not only seemingly hard, it is even known to be impossible under relatively weak assumptions. See Section~\ref{s.manygaps} for details.

Once infinitely many gaps are present, $\cS$ has infinitely many connected components. Must they be intervals or is it ever possible that the gaps of $\cS$ are dense, either in all of $\R$ or some portion of it? In other words, can $\cS$ (or a portion thereof) ever be a (generalized) Cantor set in the sense that it has empty interior?\footnote{Since all operators that we consider are unbounded, the corresponding spectra are also unbounded and hence never compact. To avoid constantly saying ``generalized Cantor set'', we will abuse terminology somewhat and say a Cantor set is a closed set with empty interior and no isolated points.} This is the next question we may ask on our quest towards possible thinness phenomena:

\begin{question}\label{q.cantor}
Can the gaps of $\cS$ be dense? In other words, is it possible that $\cS$ has empty interior?
\end{question}

The occurrence of Cantor spectra has been a hot topic since the 1980s. To exhibit this phenomenon one has to leave the class of periodic potentials, since Cantor structures are easily seen to be impossible in the periodic case. But simply taking the closure of the class of periodic potentials in a suitable way does the trick! In other words, potentials that are \emph{almost-periodic} are good candidates when looking for Schr\"odinger operators with Cantor spectra. Much of the existing theory requires $d = 1$. This is not surprising as we mentioned above that gaps at high energies are atypical for periodic potentials when $d > 1$ and approximation of almost-periodic potentials by periodic potentials is a natural tool that should also inform one's expectations about the limit spectra that arise. However, very recently it has been shown how Cantor spectra can be shown to arise even in the higher-dimensional setting. We will discuss all of these points in Section~\ref{s.cantor}.

For nowhere dense sets, which are quite thin in a topological sense, one can explore further thinness features by considering measure theoretic properties. Starting with Lebesgue measure, one may ask the following:

\begin{question}\label{q.zeroLeb}
Can the Lebesgue measure of $\cS$ vanish?
\end{question}

While the Cantor set just about anyone is first exposed to is the middle third Cantor set, which is easily seen to have zero Lebesgue measure, further explorations then show that there are also ``fat Cantor sets'' of positive Lebesgue measure. In the context of Schr\"odinger spectra the earliest examples turned out to have positive Lebesgue measure, but further, more quantitative, work then was able to give a positive answer to the previous question; see Section~\ref{s.zeroLeb} for details.

Once the Lebesgue measure of $\cS$ is shown to be zero, this naturally prompts us to consider fractal dimensions. Specifically, since Lebesgue measure on $\R$ is just the $\alpha$-dimensional Hausdorff measure for $\alpha = 1$, it is natural to lower the value of $\alpha$ until $\cS$ no longer has zero weight with respect to $h^\alpha$. This transition value if called the \emph{Hausdorff dimension} of $\cS$. There are other ways of studying whether $\cS$ has fractional dimension in a suitable sense, for example by considering the (lower or upper) box counting (or Minkowski) dimension. If a fractional value occurs for one of these dimensions, it is then of interest to look at the others and see whether they yield the same value, or whether there are essential spectra for which the different fractal dimensions indeed differ. Since the spectrum is a subset of $\R$, all these fractal dimensions of $\cS$ will take values in $[0,1]$. Once one confirms that values below $1$ are possible, one naturally asks whether dimensions all the way down to, and including, zero are possible:

\begin{question}\label{q.zerodim}
Can $\cS$ be zero-dimensional?
\end{question}

Section~\ref{s.zerodim} explains that this is indeed possible, and this will demonstrate that the existing results have pushed the possible thinness of essential spectra of Schr\"odinger operators to the extreme, in the sense of topology, measure theory, and fractal geometry. But there is an important additional aspect in pursuing this quest for extreme results. Only once zero-dimensionality has been established in the case $d = 1$, can one then address Questions~\ref{q.cantor}--\ref{q.zerodim} in the case $d > 1$. See Section~\ref{s.zerodim} for this discussion.

After all of these results, one may wonder whether there is any lower bound that one can impose on the essential spectrum.

\begin{question}\label{q.lowerbound}
Is there any sense in which the size of the set $\cS$ can be bounded from below?
\end{question}

We will take this up in Section~\ref{s.lowerbound}. Let us note that if $V(x)$ is unbounded, then it is possible for $H$ to have compact resolvent and hence for $\sigma(H)$ to consist of isolated eigenvalues of finite multiplicity which only accumulate at $+\infty$. However, in this case $\cS = \emptyset$, so we do not view this as an interesting counterexample. One may observe this explicitly for (say) the quantum harmonic oscillator ($V(x) = x^2$), which can be explicitly diagonalized in terms of Hermite functions; compare \cite[p.~142]{ReedSimon1} or \cite[Section~6.4]{Sim2015:CCA1}. Consequently, it is mainly of interest to exhibit the desired spectral results within (say) the class of bounded potentials, and we shall restrict attention to that case in this note.
\medskip

For this note, we focus primarily on operators as in \eqref{e.schrodingeroperator}, but it will sometimes be useful to discuss the related \emph{discrete Schr\"odinger operators} given by choosing a bounded potential $v:\Z^d \to \R$ and defining $h=h_v = -\Delta_{\Z^d} +v$ in $\ell^2(\Z^d)$, where
\[[\Delta_{\Z^d}\psi](n) = \sum_{\|m - n\|_1 = 1} \psi(m), \]
and
$[v\psi](n) = v(n)\psi(n)$ for $n \in \Z^d$.

\section{Can $\cS$ have gaps?}\label{s.disconnected}

As discussed, if $V \equiv 0$, the spectrum of $H$ is the half-line $[0,\infty)$, and if $V$ enjoys suitable decay at $\pm \infty$, the essential spectrum remains undisturbed. Thus, it becomes necessary to work with nonzero potentials that do not decay. Of course, if $V(x) \equiv \lambda_0$ a constant, then $\sigma(H)= [\lambda_0,\infty)$, so it is necessary to work with nonconstant potentials as well. The simplest potentials to work with are those that are periodic; we say $V: \R^d \to \R$ is \emph{periodic} if there is a basis $\{x_1,\ldots,x_d\}$ of $\R^d$ such that $V(x+x_j) \equiv V(x)$ for every $1 \le j \le d$.

It is an amazing result from the inverse theory that, in one space dimension, \emph{every} non-constant periodic potential has a gap in its spectrum.

\begin{theorem}
If $V:\R \to \R$ is periodic and nonconstant, $\cS$ has a gap.
\end{theorem}

This result is originally due to Borg~\cite{Borg1946}, with a short self-contained proof due to Hochstadt~\cite{Hochstadt1965}; see also the Historical Remarks and references in \cite{Hochstadt1965} for more.

Let us say a little bit more about the determination of the spectrum when $V$ is periodic. We emphasize that the approach we discuss here is strictly one-dimensional; we will discuss the higher-dimensional setting later. For simplicity, take $V$ to have period one, and let $u_1, u_2$ solve the initial value problem $-u'' + Vu = Eu$ subject to initial conditions
\[ \begin{bmatrix} u_1'(0) & u_2'(0) \\ u_1(0) & u_2(0) \end{bmatrix} = \begin{bmatrix} 1 & 0 \\ 0 & 1 \end{bmatrix}.\]
Then, the \emph{monodromy matrix}
\[\Phi(E) = \begin{bmatrix} u_1'(1) & u_2'(1) \\ u_1(1) & u_2(1) \end{bmatrix}\]
and its trace $D(E) := \tr [\Phi(E)]$ tell essentially the whole story about general solutions to
\begin{equation} \label{eq:disconnected:timeindepschrod}
-u''+Vu=Eu,
\end{equation} in that one has
\[ \begin{bmatrix} u'(n) \\ u(n) \end{bmatrix} = \Phi(E)^n \begin{bmatrix}u'(0) \\ u(0) \end{bmatrix}\]
for all $n \in \Z$. Thus, asymptotic characteristics of $u(x)$ are entirely dictated by behavior of large powers of $\Phi(E)$. In particular, one has the following possibilities:
\begin{itemize}
\item $D(E) \in (-2,2)$ \emph{or} $\Phi(E) = \pm I$. Every solution of \eqref{eq:disconnected:timeindepschrod} is uniformly bounded, so $E$ is a generalized eigenvalue of $H_V$, and hence $E$ belongs to the spectrum of $H_V$ by Sch'nol's Theorem \cite{Schnol1957, Simon1981JFA}.
\medskip

\item $D(E) = \pm 2$ and $\Phi(E) \not= \pm I$. Solutions of \eqref{eq:disconnected:timeindepschrod} are linearly bounded, and $E$ again belongs to the spectrum of $H_V$ by Sch'nol's theorem.
\medskip

\item $D(E) \notin [-2,2]$. There are linearly independent solutions $u_\pm$ of \eqref{eq:disconnected:timeindepschrod} such that $u_\pm$ decays exponentially at $\pm \infty$ and grows exponentially at $\mp \infty$, and $E \notin \sigma(H)$. One can appeal to Sch'nol again or use $u_\pm$ to explicitly construct the resolvent of $H$ at $E$.
\end{itemize}

\begin{figure}
\begin{center}
\begin{tikzpicture}
\draw [->,line width=.06cm] (-.5,0) -- (10.5,0);
\draw [->,line width = .06cm] (2.5,-2.5) -- (2.5,3);
\draw [smooth,samples=100,line width=.03cm,black,domain=0:3] plot (\x, {(1.49+0.01*(\x-5.9)^2)*cos((sqrt(15*\x)) r)});
\draw [smooth,samples=100,line width=.03cm,black,domain=3:10] plot (\x, {(1.49+0.01*(\x-5.9)^2)*cos((sqrt(15*\x)) r)});
\draw [-,line width = .04cm,red] (-0.5,1.5) -- (10,1.5);
\draw [-,line width = .04cm,red] (-0.5,-1.5) -- (10,-1.5);
\draw [dashed,line width =.01cm,blue] (0,1.5) -- (0,0);
\draw [dashed,line width =.01cm,blue] (.4,-1.5) -- (.4,0);
\draw [dashed,line width =.01cm,blue] (.9,-1.5) -- (.9,0);
\draw [dashed,line width =.01cm,blue] (2.28,1.5) -- (2.28,0);
\draw [dashed,line width =.01cm,blue] (2.95,1.5) -- (2.95,0);
\draw [dashed,line width =.01cm,blue] (5.9,-1.5) -- (5.9,0);
\draw [dashed,line width =.01cm,blue] (9.8,1.5) -- (9.8,0);
\draw [-,line width = .15cm,green] (0,0) -- (.4,0);
\draw [-,line width = .15cm,green] (0.9,0) -- (2.28,0);
\draw [-,line width = .15cm,green] (2.95,0) -- (9.8,0);
\node [right] at (10.5,0) {$E$};
\node [above] at (2.5,3) {$D(E)$};
\node [above] at (5.5,2) {\textcolor{red}{$D(E) = \pm 2$}};
\node [below] at (-0.2,0) {\textcolor{blue}{$a_1$}};
\node [above] at (0.4,0) {\textcolor{blue}{$b_1$}};
\node [above] at (0.9,0) {\textcolor{blue}{$a_2$}};
\node [below] at (2.28,0) {\textcolor{blue}{$b_2$}};
\node [below] at (2.95,0) {\textcolor{blue}{$a_3$}};
\node [above] at (5.9,0) {\textcolor{blue}{$b_3 = a_4$}};
\node [below] at (9.8,0) {\textcolor{blue}{$b_4$}};
\end{tikzpicture}
\end{center}
\caption{Determination  of $\cS$ from $D$. \label{fig:discriminant}}
\end{figure}
One can say a bit more about $D$. Since the spectrum of $H_V$ must be real, any value of $E$ with $D(E) \in [-2,2]$ must be real. In particular, $D$ has all real roots. Moreover, if $D'(E) = 0$ for some $E \in \R$, one must have $D(E) \notin (-2,2)$, for, if not, then one could produce nearby nonreal $E$ with $D(E) \in (-2,2)$, contradicting self-adjointness.
Thus, $\cS$ is a union of nondegenerate closed bands,
\[\cS = \bigcup_{k=1}^\infty [a_k,b_k], \]
where each band $B_k = [a_k,b_k]$ is obtained as the closure of a connected component of the preimage of $(-2,2)$ under $D$. See Figure~\ref{fig:discriminant} for a visualization of $D(E)$ and how it determines $\cS$.

In view of this picture of $D$, one sees that a gap may only close at energy $E_0$ if $D(E)$ is tangent to one of the horizontal lines at height $\pm 2$ at $E_0$. Consequently, it seems like the presence of a closed gap is an unstable phenomenon, and one might expect that a ``typical'' potential breaks all of these gaps open. We discuss this further in the next section.

Of course, the approach outlined above does not work for $d \geq 2$, but a closer inspection suggests the correct generalization. Namely, given $\theta \in \T: = \R/\Z$, we may have $D(E) = 2\cos(2\pi \theta) \in [-2,2]$ if and only if there is a solution $u$ to \eqref{eq:disconnected:timeindepschrod} such that $u(x+1) = e^{2\pi i\theta}u(x)$, and it is this notion which makes perfect sense in higher dimensions. Namely, if $V$ is periodic and $\{x_1,\ldots,x_d\}$ is a basis of $\R^d$ for which $V(x) \equiv V(x+x_j)$ for every $j$, then one seeks solutions of
\begin{equation} \label{eq:disconnected:timeindepschrodmuldtid}
-\Delta u + Vu = Eu
\end{equation}
such that
\begin{equation} \label{eq:disconnected:blochbcs} u(x+x_j) = e^{2\pi i\theta_j} u(x)\end{equation}
 for some $\theta = (\theta_1,\ldots,\theta_d) \in \T^d : = \R^d/\Z^d$. For every $\theta$, there is a countable set $E_1(\theta) \leq E_2(\theta) \leq \cdots$ such that \eqref{eq:disconnected:timeindepschrodmuldtid} enjoys a nontrivial solution for $E = E_j(\theta)$ satisfying \eqref{eq:disconnected:blochbcs} (we enumerate the $E_j$'s according to multiplicity, which is why we write ``$\leq$''). Then, the $k$th band is given by $B_k = \{E_k(\theta) : \theta \in \T^d\}$, and one has
 \[\cS = \bigcup B_k.
 \]
Each $E_k$ can be shown to be non-constant, so each $B_k$ is a non-degenerate closed interval. The primary difference between the lower- and higher-dimensional cases is that the bands could only touch in $d=1$, but when $d \geq 2$, the bands may
\emph{overlap} (i.e. one may have $b_k > a_{k+1}$). This will be responsible for some of the differences one sees in the higher-dimensional setting, including the added difficulty of constructing examples having Cantor spectrum.

\section{Can $\cS$ have infinitely many gaps?}\label{s.manygaps}

The answer to this question (in the periodic setting) depends entirely on the dimension. If $d = 1$ and, for example, $V(x) = \cos x$, then it is known that $\cS$ is a disjoint union of infinitely many compact intervals \cite[pp.~298--299]{ReedSimon4}. In other words, the spectrum has infinitely many gaps at arbitrarily high energies. This topological structure is known to be generic\footnote{Here, and throughout this note, a property that may hold at some point in a topological space is said to be generic if the set on which it holds is a dense $G_\delta$.} among all one-dimensional periodic potentials. For the sake of notation, let $C^r_{\rm per}(\R)$, $r \in \Z_+ \cup \{\infty\}$ denote those $V \in C^r(\R)$ such that $V(x+1) \equiv V(x)$, equipped with the uniform topology for $r=0$ and the usual Fr\'echet topology when $r >0$.

\begin{theorem}[Simon, 1976 \cite{S76}] \label{t:S76}
For generic $V \in C^r_{\rm per}(\R)$, $b_k < a_{k+1}$ for all $k$. In particular, $\cS$ has infinitely many open gaps.
\end{theorem}

On the other hand, the picture in higher dimensions is starkly different. In general, the arguments become quite complex, but one can get a simple idea for how things change in the following simple scenario. First, consider a given $V\in C^0_{\rm per}(\R)$. Denoting the $k$th band and gap by $B_k = [a_k,b_k]$ and $G_k=(b_k,a_{k+1})$ respectively, and using $|\cdot|$ to denote the length of an interval, we have:
\begin{align} \label{eq:infgaps:bandbound}
\sum_k \left( \pi^2(2k-1)  - |B_k|\right)^2 & < \infty \\
\label{eq:infgaps:gapbound}
\sum_k|G_k|^2 &  <\infty;
\end{align}
see \cite{GT84, GT87}. (In fact, the summands of \eqref{eq:infgaps:bandbound} and \eqref{eq:infgaps:gapbound} are very easily seen to be bounded by $4\|V\|_\infty^2$ by general perturbation theory,  and this already suffices for the purposes that we shall undertake presently). Now, let us see what insights this yields in higher dimensions. A potential $V:\R^d \to \R$ is called \emph{separable} if there are potentials $W_j : \R \to \R$ for $1 \le j \le d$ so that \[V(x) = \sum_{j=1}^d W_j(x_j).\]
Writing $\cS = \sigma(H_V)$ and $\cS_j = \sigma(H_{W_j})$, the general theory (see, for instance,  \cite[Chapter~VIII]{ReedSimon1}) shows that $\cS$ decomposes as a Minkowski sum:
\begin{equation} \label{eq:infgaps:sumspec} \cS= \sum_{j=1}^d \cS_j= \set{\sum_{j=1}^d E_j :E_j \in \cS_j}. \end{equation}
Combining the observations \eqref{eq:infgaps:bandbound}, \eqref{eq:infgaps:gapbound}, and \eqref{eq:infgaps:sumspec} shows the following: if $V:\R^d \to \R$ is a continuous separable periodic potential, then $\sigma(H_V)$  has only finitely many gaps. On the basis of this calculation, the following result was quite long-expected:

\begin{conj}[Bethe--Sommerfeld Conjecture]
Suppose $d \geq 2$. If $V:\R^d \to \R$ is periodic and is such that $H_V$ is self-adjoint, then $\sigma(H_V)$ has finitely many gaps.
\end{conj}

This result is quite difficult to prove in the non-separable case, and has only recently been concluded.
An incomplete list of contributions includes \cite{HelMoh98, Karp97, Skr79, Skr84, Skr85, Vel88}, with the current best available result in \cite{Parn2008AHP}.
Versions for discrete Schr\"odinger operators may be found in \cite{EF2017, KrugPreprint}, with the most general discrete result in \cite{HanJit}.
For more on the history and development of this problem (and the spectral theory of periodic operators writ large), we recommend Kuchment's survey \cite{Kuc2016}.

\section{Can the gaps of $\cS$ be dense?}\label{s.cantor}

From the previous discussion, the answer must be a resounding ``no'' if one restricts to the case of periodic potentials. Another scenario that was studied intensively at first was the case of random potentials. Under relatively mild assumptions, it is known that the spectrum of a random operator is the same for almost every realization and that this almost-sure spectrum is given by the closure of the union of the spectra of corresponding periodic realizations; compare \cite{KirschMart92}. In particular, if there are any gaps in the (almost-sure) spectrum, they do not accumulate. Thus, in the random and periodic cases, the associated operators (almost) always have spectra in which the gaps cannot be dense. Let us note in passing that the study of random operators is itself a deep field with a vast literature; we refer the reader to \cite{AizWar, Stollmann2001} and references therein.

However, in the 1970s and 1980s, interest began to grow in aperiodic potentials with suitable recurrence properties, motivated by several developments in math and physics. First, the study of electrons subjected to external magnetic fields (which dates back to the early part of the 20th century, compare \cite{Harper1955, Luttinger1951, Peierls1933}) led to interest in quasi-periodic models. Famously, plots of the spectra corresponding to Harper's model led to a picture now known as Hofstadter's butterfly \cite{Hofstadter1976}, one of the first fractals observed in physical scenarios. These models led to interest in the almost-Mathieu operator (AMO), which is a quasi-periodic discrete Schr\"odinger operator $h_v$ with
\begin{equation} \label{eq:cantor:amodef} v(n) = v_{\lambda,\alpha,\omega}(n)= 2\lambda \cos(2\pi(n\alpha+\omega))\end{equation}
for $n \in \Z$, where $\lambda,\alpha, \omega \in \R$ are called the \emph{coupling constant}, \emph{frequency}, and \emph{phase}, respectively. We recommend \cite{D2017, MarxJito2017} for more about such quasi-periodic operators, their history, and connections to physics. Of course, if $\alpha$ is rational or $\lambda = 0$, the resulting $v$ is periodic and hence can be handled by a formalism analogous to that described in Section~\ref{s.disconnected}, so the most interesting case is that in which $\alpha$ is irrational and $\lambda \not=0$. The second major physical development that spurred mathematical study of aperiodic operators was the discovery of quasicrystals by Shechtman \textit{et al.}\ in the early 1980s \cite{SBGC1984PRL}. We will not give a comprehensive account of mathematical quasicrystals, instead referring the reader to \cite{BaaGri1} and references therein for background. The central operator-theoretic model of a quasicrystal is the \emph{Fibonacci Hamiltonian}, a discrete Schr\"odinger operator with potential given by
\begin{equation}
\label{eq:cantor:FHdef}v(n) = v_{\lambda,\omega}(n)= \lambda \chi_{[1-\alpha,1)}(n\alpha + \omega \ \mathrm{mod} \ 1),
\end{equation}
where $\lambda>0$ and $\omega \in \R$ are parameters, again called the coupling constant and the phase. Here, $\alpha$ is fixed to be $\frac{1}{2} (\sqrt{5}-1)$, the inverse of the golden mean (however, one can also choose other irrational values of $\alpha$, leading to \emph{Sturmian} potentials).

Let us briefly comment on the similarities and differences between \eqref{eq:cantor:amodef} and \eqref{eq:cantor:FHdef}. Both are aperiodic, but both are also ``not too far'' from being periodic, however in vastly different ways. As an almost-periodic potential, \eqref{eq:cantor:amodef} exhibits $\varepsilon$-almost periods, that is, nontrivial translations $t$ such that $\|v - v(\cdot -t)\|_\infty < \varepsilon$. The Fibonacci potential \eqref{eq:cantor:FHdef} is not almost-periodic, indeed has no nontrivial $\varepsilon$-almost-period for any $\varepsilon \leq \lambda$. Instead, there are nontrivial translations $t_1, t_2, \ldots$ and $R_j \to \infty$ such that $v$ and $v(\cdot-t_j)$ agree on $[-R_j,R_j]$. To summarize, \eqref{eq:cantor:amodef} enjoys imprecise repetitions at a global scale, while \eqref{eq:cantor:FHdef} enjoys exact local repetitions. These differences inform the different properties of these operators, as well as the different tools and techniques that have been brought to bear in their study.

To study an operator with a potential such as \eqref{eq:cantor:amodef} (or a suitable analogue in $L^2(\R)$), a natural tactic is to approximate $\alpha$ by rationals, study the resulting periodic operators, and try to prove suitable approximation results. The problem is that such approximations are inherently \emph{local}. In particular, these approximations will be close to the original aperiodic potential on a finite window, and will not at all be good approximations outside that window. To put this into the language of operator theory, the resulting periodic approximants approximate the aperiodic operator \emph{strongly} but not in \emph{norm}, so one is naturally interested in what operators enjoy a periodic approximation in the norm topology. To wit: a potential $V \in C(\R^d)$ is called (uniformly) \emph{limit-periodic} if it may be uniformly approximated by periodic elements of $C(\R^d)$. For the case of limit-periodic operators, one might expect their spectral properties to correspond closely to those in the periodic case, but we will see that they can diverge quite far. Let $\LP(\R^d)$ denote the set of limit-periodic $V \in C(\R^d)$, which is a complete metric space in the uniform metric (note however that it is not closed under sums and hence is not a Banach space). The first result on Cantor spectrum in the limit-periodic case is due to Moser:

\begin{theorem}[Moser, 1981 \cite{Moser1981}] \label{t:Moser81}
There exist $V \in \LP(\R)$ such that $\cS$ is a Cantor set.
\end{theorem}

Very shortly after Moser's work, Avron--Simon \cite{AS81} showed that Cantor spectrum is generic in $\LP(\R)$. Let us also mention related work of Chulaevskii \cite{Chul81, Chul84} and Pastur--Tkachenko \cite{PT1, PT2}. For more about limit-periodic operators in general, see the survey \cite{DF:LPSurv} and references therein.

Once one understands how the process of periodic approximation works, one is no longer surprised by theorems showing that Cantor spectrum is generic in the limit-periodic setting. There are a few basic heuristics: First, one can show from the definitions that the spectrum never has isolated points. Next, one can show that Cantor spectrum is a $G_\delta$ phenomenon, which is essentially an exercise in point-set topology and perturbation theory. Concretely, if $I$ is an open interval, the set of $V$ for which $I \cap \rho(H_V) \not=\emptyset$ is easily shown to be open. Thus, the intersection of all such sets of potentials as $I$ ranges over intervals having rational endpoints is a $G_\delta$ and is rather clearly the set of potentials enjoying nowhere-dense spectrum. Finally, one needs to show that Cantor spectrum is dense in the space $\LP(\R)$. This is more technical, but the outline is the following. Start with a periodic potential, and successively produce periodic perturbations of higher and higher periods for which all gaps that can open do open (by using Theorem~\ref{t:S76}), and in such a way that gaps opened in previous steps do not close, taking care that the gaps also cannot close in the limit.

Let us observe that the outline suggested in the previous paragraph does \emph{not} work for dimensions $d>1$. Concretely, Theorem~\ref{t:S76} only holds for $d=1$, and fails spectacularly for $d>1$, in view of the Bethe--Sommerfeld conjecture. In particular, the presence of Cantor spectrum in higher dimensions is necessarily more difficult to prove than in the 1D case. Nevertheless, we will see in Section~\ref{s.zerodim} that suitable quantitative refinements of the 1D Cantor results can be used to generate Cantor spectrum in higher dimensions.

Since we mentioned the AMO, let us note that it was later shown to have Cantor spectrum for all $\lambda\neq 0$, all $\omega \in \R$, and all irrational $\alpha$ \cite{AviJit2009}. Similarly, Cantor spectrum has been established for the Fibonacci (as well as the more general Sturmian) case for all non-trivial parameter values \cite{BelIocScoTes1989, Suto1989JSP}. See \cite{D2017, MarxJito2017} for further comments on the history and earlier partial results.

\section{Can the Lebesgue measure of $\cS$ vanish?}\label{s.zeroLeb}

As soon as one has a Cantor spectrum, one is interested in the measure properties. In particular, even if the gaps are dense, the spectrum may still be a ``fat'' Cantor set and hence have positive measure (or even infinite measure, since $\cS$ is unbounded).
 Returning briefly to the discrete setting, the first two operators known to have zero-measure spectrum of physical interest were the critical AMO (\eqref{eq:cantor:amodef} with $\lambda = 1$) and the Fibonacci Hamiltonian \cite{Suto1989JSP}, with potential given by \eqref{eq:cantor:FHdef}.
On one hand, the zero-measure result for the AMO is very sensitive to the coupling constant, since the measure of the spectrum (for irrational $\alpha$) is known to be given by $|4-4|\lambda||$ for every $\lambda$ (see e.g.\ \cite{AviKri2006, JitoKras2002, Last1993, Last1994}).
In particular, as soon as one perturbs the coupling constant, the zero-measure property is destroyed, even though the Cantor property persists.
By way of contrast, the Fibonacci Hamiltonian exhibits zero-measure Cantor spectrum for \emph{arbitrary} $\lambda \neq 0$, and hence, the zero-measure phenomenon is more robust in this scenario.

This discussion of stability is significant from the point of view of the desired results for operators of the form $-\Delta+V$, which are unbounded. In particular, the high-energy region for $H_V$ is analogous to the regime of small coupling constant in the discrete case. In particular, one should not expect a zero-measure result to hold for an analogue of the AMO in $L^2(\R)$, so, in order to generate examples of the desired phenomenon for Schr\"odinger operators in $L^2(\R)$, one should look in (say) the Fibonacci setting, where the desired phenomenon is robust in the small coupling region.

To spare the reader technicalities and notation, let us work in a setting that is simpler than the most general possible setup, but which still exhibits the behaviors that we want to discuss. Let $u$ denote the $0$-$1$ Fibonacci sequence, given as:
\[u(n) = u(n,\omega) =  \chi_{[1-\alpha,1)}(n\alpha + \omega \ \mathrm{mod} \ 1)\]
with $\alpha = \frac{1}{2}(\sqrt{5}-1)$ as before.
Choosing $f_0 \neq f_1 $ in $L^2([0,1))$, one can generate a Fibonacci Schr\"odinger operator by replacing each $0$ in $u$ by $f_0$ and each $1$ by $f_1$ and attaching a coupling constant\footnote{Of course, one could simply absorb $\lambda$ into the $f_j$, but it is of interest to fix $f_j$ first and then to observe the scaling behavior as the coupling constant is varied.} $\lambda$; more precisely, take:
\begin{equation} \label{eq:contFHdef} V_{\omega,\lambda}(n) = \lambda \sum_{n \in \Z} f_{u(n,\omega)}(\cdot - n) \chi_{[n,n+1)}. \end{equation}

\begin{theorem}[D.--F.--Gorodetski, 2014 \cite{DFG2014}]
For every $f_0 \neq f_1$, every $\lambda \neq 0$, and every $\omega \in \R$, $\cS$ is a Cantor set of zero Lebesgue measure.
\end{theorem}

The result of \cite{DFG2014} is more general than this: one can replace $\alpha$ by an arbitrary irrational number; in fact, one can replace these aperiodic $0$-$1$ sequences by any element of a minimal subshift satisfying Boshernitzan's criterion for unique ergodicity \cite{Bosh1984, Bosh1992}, subject to a suitable local recognizability criterion on the local potential pieces; see \cite{DFG2014} for details. See also \cite{LSS} for similar results that were proved independently. We specialize to the Fibonacci case because of the physical motivation and because our later results about dimension will require it.

Additionally, the zero-measure property also holds for the case of limit-periodic potentials, and remains stable in the coupling constant for generic potentials, first shown in the discrete case by Avila~\cite{avila2009CMP}.

\begin{theorem}[D.--F.--Lukic, 2017 \cite{DFL2017}]
For generic $V \in \LP(\R)$, $\sigma(H_{\lambda V})$ is a zero-measure Cantor set for every $\lambda \neq 0$.
\end{theorem}

The latter claim is somewhat surprising, since one knows that $\sigma(H_{\lambda V})$ converges to $\sigma(H_0) = [0,\infty)$ in the Hausdorff metric as $\lambda \to 0$. In particular, as $\lambda \downarrow 0$, the resolvent set of $H_{\lambda V}$ remains a dense, full-measure subset, even though the length of every gap goes to zero at least as fast as $2\lambda \|V\|_\infty$.

\section{Can $\cS$ be zero-dimensional?}\label{s.zerodim}

Drawing on ideas from Avila \cite{avila2009CMP}, one can refine the Cantor spectrum arguments in a quantitative way. The first step is to go from qualitative behavior to quantitative behavior in the periodic case. Concretely, we know (from Theorem~\ref{t:S76}) that generic periodic potentials open all possible gaps in dimension one. One of the key insights of Avila in the discrete setting was a mechanism that one can use to generate not only gaps, but to show that the spectrum of an operator of period $p$ can be made exponentially small (i.e. $\lesssim e^{-cp}$ for a suitable constant $c>0$). Altering the original construction of limit-periodic operators to generate periodic approximants with such measure estimates at each stage, one can control the fractal dimensions of the limiting spectrum. Of course, in the continuum case, it does not make sense to  speak of measure bounds on the whole spectrum of a periodic operator (which necessarily has infinite measure, compare \eqref{eq:infgaps:bandbound}), so some care is needed. Nevertheless, one can show

\begin{theorem}[D.--F.--Lukic, 2017 \cite{DFL2017}, D.--F.--Gorodetski, 2019 \cite{DFG2019}] \label{t:dfgl}
Within $\LP(\R)$, the set of $V$ for which $\dim_{\rm H}(\sigma(H_V))=0$ is dense. In fact, for a dense set of $V \in \LP(\R)$, one has $\dim_{\rm H}(\sigma(H_{\lambda V})) = \dim_{\rm B}^-(\sigma(H_{\lambda V})) = 0$ for every $\lambda \neq 0$.
\end{theorem}

The result about Hausdorff dimension is from \cite{DFL2017}, and the result abou the lower box-counting dimension is from \cite{DFG2019}. The strengthening from Hausdorff to box-counting dimension allows us to conclude answers to all of our questions about $\cS$ in higher dimensions.

\begin{theorem}[D.--F.--Gorodetski, 2019 \cite{DFG2019}] \label{t:dfg2019}
There exist $V \in \LP(\R^d)$ for any $d \geq 2$ such that $\cS$ is a Cantor set having zero Lebesgue measure, zero Hausdorff dimension, and zero lower box-counting dimension.
\end{theorem}

The key fact that enables one to make the transition is the following observation. If $A$ is a compact set, then
\begin{equation} \label{eq:dimsum} \dim_{\rm B}^- \left(\sum_{j=1}^d A\right) \leq d\dim_{\rm B}^-(A).\end{equation}
In particular, if $A$ has zero lower box-counting dimension, so too does $\sum_{j=1}^d A$ for any $d$. Thus, Theorem~\ref{t:dfg2019} is proved by using Theorem~\ref{t:dfgl} to construct a 1D potential $W$ for which the associated spectrum has lower box-counting dimension zero, constructing a $d$-dimensional separable potential of the form
\[V(x) = \sum_{j=1}^d W(x_j),\]
and then invoking \eqref{eq:dimsum} and \eqref{eq:infgaps:sumspec}.
We are sweeping a few niceties about fractal dimensions of sums of unbounded sets under the rug here, but they are not hard to deal with, see, e.g., the appendix to \cite{DFG2019}.

Let us note that one really does need the limit-periodic setting here.
One can show that the critical AMO exhibits zero Hausdorff dimension for some frequencies \cite{LastShamis2016}; however, we remind the reader that because these only hold for $\lambda = 1$, one cannot expect to be able to use them to deduce consequences for Schr\"odinger operators in $L^2(\R)$.
Moreover, one can prove dimensional results for the Fibonacci Hamiltonian, but one cannot zero out its Hausdorff dimension.
In particular, one has the following for potentials $V_{\omega,\lambda}$ as in \eqref{eq:contFHdef}.
Let us denote $\cS_\lambda$ the essential spectrum corresponding to $V_{\omega,\lambda}$ (the omission of $\omega$ from the notation reflects the fact that $\cS$ does not depend on $\omega$, which can be seen from strong approximation arguments).

\begin{theorem}[D.--F.--Gorodetski, 2014 \cite{DFG2014}, F.--Mei, 2018 \cite{FM2018}]\label{t:dfgm}
For any $f_0 \neq f_1$ in $L^2([0,1))$ and any $\lambda>0$, one has
\[\dim_{\rm H}^{\rm loc}(\cS_\lambda,E) >0\]
for every $E \in \cS_\lambda$. In the small-coupling limit, one has:
\[\lim_{\lambda \to 0 } \inf_{E \in \cS_\lambda} \dim_{\rm H}^{\rm loc} (\cS_\lambda,E) = 1.\]
Moreover, for any $\lambda>0$, in the high-energy limit, we have
\[ \lim_{E_0 \to \infty} \inf_{E \in \cS_\lambda \cap [E_0,\infty)} \dim_{\rm H}^{\rm loc} (\cS_\lambda,E) = 1.\]
In particular, the global Hausdorff dimension of $\cS_\lambda$ is always equal to $1$.
\end{theorem}

The asymptotic statements for the dimension were proved for $f_0 \equiv 0$ and $f_1 \equiv 1$ in \cite{DFG2014} and extended to general $f_j$ in \cite{FM2018}.
Let us now move to the higher dimensional setting. We fix $f_0 \equiv 0$, $f_1 \equiv \lambda$, and denote by $\cS_\lambda$ the spectrum of $H_{V_{\omega,\lambda}}$ as before. For $d \geq 2$, write \[\cS_\lambda^{(d)} = \sum_{j=1}^d \cS_\lambda,\]
the $d$-fold sum (cf.\ \eqref{eq:infgaps:sumspec}). In particular, $\cS_\lambda^{(d)}$ is the spectrum of $H_{\lambda,\omega}^
{(d)} := -\Delta + V_{\omega,\lambda}^{(d)}$, where
\[V_{\omega,\lambda}^{(d)}(x) = \sum_{j=1}^d V_{\omega,\lambda}(x_j), \quad x \in \R^d.\]
The dimensional results in Theorem~\ref{t:dfgm} suggest (but do not imply!) that one has the following picture for the multidimensional operator $H_{\omega,\lambda}^{(d)}$: the spectrum undergoes a transition from Cantor-like at low energies to containing intervals and perhaps eventually a half-line for large energies. This was recently proved.

\begin{theorem}[D.--F.--Gorodetski, 2020+ \cite{DFG2020}]
Let $d\geq 2$. For any $\lambda>0$, $\cS_\lambda^{(d)}$ contains a half-line. There exists $C(d)>0$ such that, for $\lambda \geq C(d)$, $\cS_\lambda^{(d)}$ enjoys Cantor structure near the ground state, that is, there is an $E_0$ such that $(-\infty,E_0] \cap \cS_\lambda^{(d)}$ is a nonempty Cantor set.
\end{theorem}

\section{When Does it End?} \label{s.lowerbound}

With these increasingly fine results that show that the spectrum can be made arbitrarily small in the sense of Lebesgue measure and fractal dimensions, one wonders whether there is any hard limit at all to how far one can push this. In other words, is there some notion of the size of a set for which one can prove a \emph{lower} bound on the size of the spectrum, valid for a large class of (bounded) potentials? Of course, as we already mentioned one does need to insert the word ``bounded'' here, since there are operators with unbounded potentials having compact resolvent, for which $\cS = \emptyset$. It turns out that potential theory in the complex plane offers a partial answer. Let us emphasize that the results in this section only apply to one-dimensional operators. At this time, we are unaware of any similar results for higher-dimensional operators.

Let us briefly recall the definitions and refer the reader to \cite{Ransford, StaTot} for textbook treatments. For a compact set $K \subseteq \C$, write $M_1(K)$ for the set of all Borel probability measures supported on $K$. For $\mu \in M_1(K)$, the (logarithmic) energy of $\mu$ is given by
\[\mathcal{E}(\mu) = - \iint \log|x-y| \, d\mu(x) \, d\mu(y),\]
with $\mathcal{E}(\mu) = +\infty$ an allowed value. The (logarithmic) capacity of $K$ is
\[ \mathrm{Cap}(K) := \sup\{ e^{-\mathcal{E}(\mu)} : \mu \in M_1(K)\}, \]
where $e^{-\infty} = 0$ by convention. A set with $\mathrm{Cap}(K) = 0$ is known as a \emph{polar set}; such sets play the role of ``negligible'' sets in potential theory, analogous to the role of measure-zero subsets in measure theory. Polar sets are very small in the senses that we have discussed herein: they have Hausdorff dimension zero (hence also Lebesgue measure zero) \cite{Carleson1967}; see also the exercises to \cite[Section~3.6]{Sim2015:CCA3}. In terms of these objects, there is a particularly clean lower bound for discrete Schr\"odinger operators in 1D: one must have $\mathrm{Cap}(\cS) \geq 1$.
In particular, the spectrum of such a discrete operator may have Hausdorff dimension zero, but it may never be a polar set.
This follows from a work of Szeg\H{o} \cite{Szego1924}; see also \cite{Simon2007IPI} for a modern discussion. Potential theory for unbounded sets is a good deal more delicate, so it is much harder to phrase things in such a concise fashion. Nevertheless, Eichinger--Lukic have a recent preprint that analyzes Schr\"odinger operators in $L^2(\R)$ from this point of view \cite{EichLukPreprint}.

\end{document}